\font\tenmsa=msam10
\font\tenmsb=msbm10

\font\largebf=cmbx10 scaled\magstep2
\def\all{\hbox{for all}}
\def\All{\hbox{For all}}
\def\and{\hbox{and}}
\def\bra#1#2{\langle#1,#2\rangle}
\def\Bra#1#2{\big\langle#1,#2\big\rangle}
\def\cite#1\endcite{[#1]}
\def\dom{\hbox{\rm dom}}
\def\dbs{^{**}}

\def\f#1#2{{#1 \over #2}}
\def\half{{\textstyle\f12}}
\def\infn{\inf\nolimits}
\def\Lt{{\wt L}}

\def\lin{{\rm lin}}
\def\lr{\Longrightarrow}
\def\niff{\Longleftrightarrow}
\def\qed{\hfill\hbox{\tenmsa\char03}}
\def\qlr{\quad\lr\quad}
\def\qLt{q_{\wt L}}
\def\quand{\quad\and\quad}
\def\r{\hbox{\tenmsb R}}
\def\rbar{\,]{-}\infty,\infty]}
\long\def\slant#1\endslant{{\sl#1}}
\def\st{\hbox{such that}}
\def\T{{\rm T}}
\def\wh{\widehat}
\def\wt{\widetilde}
\def\[{\big[}
\def\]{\big]}
\def\({\big(}
\def\){\big)}
\def\defSection#1{}
\def\defCorollary#1{}
\def\defDefinition#1{}
\def\defExample#1{}
\def\defLemma#1{}
\def\defNotation#1{}
\def\defProblem#1{}
\def\defRemark#1{}
\def\defTheorem#1{}
\def\locno#1{}
\def\meqno#1{\eqno(#1)}
\def\nmbr#1{}
\def\Proof{\medbreak\noindent{\bf Proof.}\enspace}
\def\Proo{{\bf Proof.}\enspace}
\def\Signoff{}
%%%%%%%%%%%%%%%%%%%%%%%%%%%%%%%%%%%%%%%%%%%%%%%%%%%%%%%%%%%%%%%%%%%%%%%
\def \INTsec{1}
\def \BSsec{2}
\def \PCdef{2.1}
\def \POLlem{2.2}
\def \SNsec{3}
\def \SNdef{3.1}
\def \IOTAone{1}
\def \DIFFQ{2}
\def \QCONTone{3}
\def \QMAXone{4}
\def \Hex{3.2}
\def \LINlem{3.3}
\def \SNDsec{4}
\def \DUALdef{4.1}
\def \DUALone{5}
\def \DUALtwo{6}
\def \LAlem{4.2}
\def \DUALthree{7}
\def \DUALfour{8}
\def \LINsec{5}
\def \ALHORlem{5.1}
\def \QCSTARlem{5.2}
\def \NORMlem{5.3}
\def \INFLrem{5.4}
\def \INFlem{5.5}
\def \INFone{9}
\def \INFtwo{10}
\def \INFfour{11}
\def \INFthree{12}
\def \INFONEcor{5.6}
\def \INFTWOcor{5.7}
\def \INFTrem{5.8}
\def \INFthm{5.9}
\def \INFCrem{5.10}
\def \INFcor{5.11}
\def \INFCone{13}
\def \INFCtwo{14}
\def \INFCthree{15}
\def \INFCfour{16}
\def \SURcor{5.12}
\def \LINMONsec{6}
\def \EEex{6.1}
\def \NIdef{6.2}
\def \RECASTlem{6.3}
\def \MONthm{6.4}
\def \MONrem{6.5}
\def \MONcor{6.6}
\def \MONone{17}
\def \MONtwo{18}
\def \MONfour{19}
\def \BBcor{6.7}
\def \ARENS{1}
\def \BB{2}
\def \BBWYLIN{3}
\def \BBWYFP{4}
\def \BBWYBB{5}
\def \BRON{6}
\def \GMS{7}
\def \KN{8}
\def \ASD{9}
\def \FENCHEL{10}
\def \RANGE{11}
\def \HBM{12}
\def \SSDMON{13}
\def \LINPOS{14}
\def \YAO{15}
\def \ZBOOK{16}
%%%%%%%%%%%%%%%%%%%%%%%%%%%%%%%%%%%%%%%%%%%%%%%%%%%%%%%%%%%%%%%%%%%%%%%
%
\magnification 1200
\headline{\ifnum\folio=1
{\hfil{\largebf Polar subspaces and automatic maximality}\hfil}
\else\centerline{\rm {\bf Polar subspaces and automatic maximality}}\fi}
\bigskip
\centerline{\bf S. Simons}
\bigskip
\centerline{\sl Department of Mathematics, University of California}
\centerline{\sl Santa Barbara, CA 93106-3080, U.S.A.}
\centerline{\sl simons@math.ucsb.edu}
\medskip
\centerline{\bf Abstract}
\medskip
\noindent
This paper is about certain linear subspaces of Banach SN spaces (that is to say Banach spaces which have a symmetric nonexpansive linear map into their dual spaces).   We apply our results to monotone linear subspaces of the product of a Banach space and its dual.   In this paper, we establish several new results and also give improved proofs of some known ones in both the general and the special contexts.
%:Introduction
\defSection \INTsec
\medbreak
\centerline{\bf \INTsec.\quad Introduction}
\medskip
\noindent
This paper is a sequel to \cite\LINPOS\endcite, in which we developed the theory of Banach SN spaces (that is to say Banach spaces which have a symmetric nonexpansive linear map into their dual spaces) and ``$L$--positive'' sets, and applied it to a number of results about monotone linear subspaces of the product of a Banach space and its dual.   Banach SN spaces, $L$--positive sets and Banach SN duals are defined formally in Sections \SNsec\ and \SNDsec.
\par 
Section \LINsec\ is about closed linear $L$--positive subspaces of those Banach SN spaces that possess Banach SN duals.   The basic technical result of this paper is Lemma \INFlem.    The novelty of our proof is that we use a penalty function of the form $2\eta\|\cdot\|$ rather than $\half\|\cdot\|^2$.   This not only produces significantly sharper results, but also leads to simpler proofs than those previously given for results of this kind.   We actually use two special cases of Lemma \INFlem:  Corollary \INFONEcor\ and Corollary \INFTWOcor.
\par
The most important consequences of Lemma \INFlem\ are Theorem \INFthm\ and Corollary \INFcor.   In Remarks \INFTrem\ and \INFCrem, we give details of the connections between these results and those already proved in \cite\LINPOS\endcite.
\par
The results on monotone subspaces in Theorem \MONthm, Corollary \MONcor\ and Corollary \BBcor\ all have one--line proofs using the corresponding results on Banach SN spaces, and these latter results have much simpler proofs that those that we have seen in the literature in the monotone case.   In Theorem \MONthm\ we establish the following result on ``automatic maximality'', which we have not seen in the literature: \slant Let $A$ be a closed linear monotone subspace of product of a real Banach space with its dual, and the adjoint subspace of $A$ be monotone. Then $A$ and its adjoint are both maximally monotone in their respective spaces\endslant.    Theorem \MONthm\ and Corollary \MONcor\ are discussed in Remark \MONrem.
\par
The motivating result behind this analysis is the classical Brezis--Browder theorem (see \cite\BB\endcite): \slant Let $A$ be a closed linear monotone subspace of the product of a reflexive real Banach space with its dual space.  Then the adjoint subspace of $A$ is monotone if, and only if, $A$ is maximally monotone\endslant.  This was generalized by Bauschke, Borwein, Wang and Yao (see \cite\BBWYLIN\endcite\ and \cite\BBWYBB\endcite): \slant If $A$ is a closed linear monotone subspace of the product of a real Banach space with its dual space then the adjoint subspace of $A$ is monotone if, and only if, $A$ is maximally monotone of type (NI)\endslant.   See Corollary \MONcor\ for our generalization of this and Corollary \BBcor\ for a proof of the original Brezis--Browder theorem.
\par
We use some basic tools of convex analysis (that is to say Rockafellar's formula for the subdifferential of a sum, and the Br\o ndsted--Rockafellar theorem).   In addition, we use two simple preliminary results, Lemma \POLlem\ and Lemma \LINlem.  
\par
The author would like to thank Liangjin Yao for sending him the preprints \cite\BBWYLIN\endcite, \cite\BBWYFP\endcite, \cite\BBWYBB\endcite\ and \cite\YAO\endcite, and Maicon Marques Alves for sending him the preprint \cite\ASD\endcite.
%:Banach space notation --- Polar subspaces
\defSection\BSsec
\medbreak
\centerline{\bf \BSsec.\quad Banach space notation --- Polar subspaces}
\medskip
\noindent
We start off by introducing some Banach space notation.
%:  Definition \PCdef
\defDefinition \PCdef
\medbreak
\noindent
{\bf Definition \PCdef.}\enspace If $X$ is a nonzero real Banach space and $f\colon\ X \to \rbar$, we say that $f$ is \slant proper\endslant\ if there exists $x \in X$ such that $f(x) \in \r$.   We write $X^*$ for the dual space of $X$ \big(with the pairing $\bra\cdot\cdot\colon X \times X^* \to \r$\big) and $X\dbs$ for the bidual of $X$ \big(with the pairing $\bra\cdot\cdot\colon X^* \times X\dbs \to \r$\big).   If $x \in X$, we write $\wh x$ for the canonical image of $x$ in $X\dbs$, that is to say\quad $x \in X\ \and\ x^* \in X^* \lr\bra{x^*}{\wh x} = \bra{x}{x^*}$.\quad We write ``lin'' for ``linar span of''.   If $Y$ is a linear subspace of $X$, we write $Y^0$ for the ``polar subspace of $Y$'', that is to say the linear subspace  $\big\{x^* \in X^*\colon\ \bra{Y}{x^*} = \{0\}\big\}$ of $X^*$.
\par
A word is in order about the proof of Lemma \POLlem\ below.   The proof we give is an adaptation of a classical proof of a result about linear functionals \(see, for instance, \cite\KN, Theorem 1.3, p,\ 7\endcite\).   However, it is worth pointing out that when $Y$ is closed the result also follows from the Robinson--Attouch--Brezis theorem applied to the indicator functions of $Y$ and $Z$.   We will use Lemma \POLlem\ in Theorem \INFthm(b).
%:  Lemma \POLlem
\defLemma \POLlem
\medbreak
\noindent
{\bf Lemma \POLlem.}\enspace\slant Let $Y$ be a linear subspace of a real Banach space $X$ and $x^* \in X^*$.   Write $Z$ for the linear subspace $\big\{x \in X\colon\ \bra{x}{x^*} = 0\big\}$ of $X$.   Then $(Y \cap Z)^0 = \lin\{x^*,Y^0\}$.\endslant
%:   Proof of Lemma \POLlem
\Proof Since $x^* \in Z^0 \subset (Y \cap Z)^0$ and $Y^0 \subset (Y \cap Z)^0$, and $(Y \cap Z)^0$ is a linear subspace, it is immediate that $\lin\{x^*,Y^0\} \subset (Y \cap Z)^0$.
\par
We now prove the opposite inclusion.   If $Y \subset Z$ then $Y \cap Z = Y$, so $(Y \cap Z)^0 = Y^0 \subset \lin\{x^*,Y^0\}$, as required.   We now come to the more interesting case, where  $Y \not\subset Z$.   Let $z^*$ be an arbitrary element of $(Y \cap Z)^0$.   We fix $x \in Y \setminus Z$, from which $\bra{x}{x^*} \ne 0$.  Let $y$ be an arbitrary element of $Y$, and define
$$z:= y - \f{\bra{y}{x^*}}{\bra{x}{x^*}}x \in X.$$
Since $y,x \in Y$, it follows that $z \in Y$.  Since
$$\bra{z}{x^*} = \bra{y}{x^*} - \f{\bra{y}{x^*}}{\bra{x}{x^*}}\bra{x}{x^*} = 0,$$
we have $z \in Z$.   Thus $z \in Y \cap Z$.  Then, since $z^* \in (Y \cap Z)^0$, 
$$0 = \bra{z}{z^*} = \bra{y}{z^*} - \f{\bra{y}{x^*}}{\bra{x}{x^*}}\bra{x}{z^*} = \bigg\langle y,z^* - \f{\bra{x}{z^*}}{\bra{x}{x^*}}x^*\bigg\rangle.$$
Since this holds for all $y \in Y$,
$$z^* - \f{\bra{x}{z^*}}{\bra{x}{x^*}}x^* \in Y^0.$$ 
But then
$$z^* \in \f{\bra{x}{z^*}}{\bra{x}{x^*}}x^* + Y^0 \in \lin\{x^*,Y^0\}.$$
This completes the proof of Lemma \POLlem.\qed
%:Banach SN spaces
\defSection\SNsec
\medbreak
\centerline{\bf \SNsec.\quad Banach SN spaces}
\medskip
\noindent
We now introduce Banach SN spaces \(these were called \slant Banach SNL spaces\endslant\ in \cite\LINPOS\endcite\) and $L$--positive sets \(\cite\LINPOS, Section 2, pp.\ 604--606\endcite\).
%:  Definition \SNdef
\defDefinition \SNdef
\medskip
\noindent
{\bf Definition \SNdef.}\enspace Let $B$ be a nonzero real Banach space.   A \slant SN map on $B$\endslant\  (``SN'' stands for ``symmetric nonexpansive''), is a linear map $L\colon\ B \to B^*$ such that
$$\|L\| \le 1\quad\and,\quad\all\ b,c\in B,\quad \bra{b}{Lc} = \bra{c}{Lb}.\meqno\IOTAone$$
A \slant Banach SN space\enspace $(B,L)$\endslant\ is a nonzero real Banach space $B$ together with a SN map $L\colon\ B \to B^*$.   We define the function $q_L\colon\ B \to \r$ by\quad $q_L(b) := \half\bra{b}{Lb}$\quad($b \in B$)\quad   (``$q$'' stands for ``quadratic'').
It follows from (\IOTAone) that, for all $a,b \in B$,
$$|q_L(a) - q_L(b)| = \half|\bra{a}{La} - \bra{b}{Lb}| =\break\half\big|\Bra{a - b}{L(a + b)}\big| \le \half\|a - b\|\|a + b\|,\meqno\DIFFQ$$
and so
$$q_L\ \hbox{is continuous on}\ B.\meqno\QCONTone$$
Now let $(B,L)$ be a Banach SN space and $A \subset B$.   We say that $A$ is \slant$L$--positive\endslant\ if $A \ne \emptyset$ and\quad $b,c \in A \lr q_L(b - c) \ge 0$.\quad We say that $A$ is \slant maximally $L$--positive\endslant\ if $A$ is $L$--positive and $A$ is not properly contained in any other $L$--positive set.   In this case,
$$b \in B \lr \infn_{a \in A} q_L(a - b) \le 0.\meqno\QMAXone$$
(If $d \in B \setminus A$ then the maximality gives us $a \in A$ such that $q_L(a - d) < 0$, while if $d \in A$ then $q_L(d - d) = 0$.)
\medskip
There are many examples of Banach SN spaces and their associated $L$--positive sets.   The following are derived from \cite\SSDMON, Examples 2.3, 2.5, pp.\ 230--231\endcite.   More examples can be derived from  \cite\SSDMON, Remark 6.7, p.\ 246\endcite\ and \cite\GMS\endcite.   The significant example which leads to results on monotonicity appeared in \cite\SSDMON, Example 6.5, p.\ 245\endcite. We will return to it here in Example \EEex.
%:  Example \Hex
\defExample \Hex 
\medbreak
\noindent
{\bf Example \Hex.}\enspace Let $B$ be a Hilbert space with inner product $(b,c) \mapsto \bra{b}{c}$ and $L\colon B \to B$ be a nonexpansive self--adjoint linear operator.   Then $(B,L)$ is a Banach SN space.   Here are three special cases of this example:
\smallbreak
(a)\enspace If, for all $b \in B$, $Lb = b$ then every subset of $B$ is $L$--positive.
\smallbreak
(b)\enspace If, for all $b \in B$, $Lb = -b$ then the $L$--positive sets are the singletons.
\smallbreak
(c)\enspace If $B = \r^3$ and $L(b_1,b_2,b_3) = (b_2,b_1,b_3)$ and $M$ is any nonempty monotone subset of $\r \times \r$ (in the obvious sense) then $M \times \r$ is a $L$--positive subset of $B$.   The set $\r(1,-1,2)$ is a $L$--positive subset of $B$ which is not contained in a set $M \times \r$ for any monotone subset of $\r \times \r$.   The helix $\big\{(\cos\theta,\sin\theta,\theta)\colon \theta \in \r\big\}$ is a $L$--positive subset of $B$, but if $0 < \lambda < 1$ then the helix $\big\{(\cos\theta,\sin\theta,\lambda\theta)\colon \theta \in \r\big\}$ is not.
\smallbreak
(d)\enspace If $B = \r^3$ and $L(b_1,b_2,b_3) = (b_2,b_3,b_1)$ then, since (\IOTAone) fails, $(B,L)$ cannot be a Banach SN space. 
\medskip
We will use the following special property of $L$--positive linear subspaces in Theorem \INFthm(b):
%:  Lemma \LINlem
\defLemma \LINlem
\medbreak
\noindent
{\bf Lemma \LINlem.}\enspace\slant Let $A$ be a linear subspace of $B$, $b \in B$ and $\{b\} \cup A$ be $L$--positive.   Then the linear subspace $\lin\big\{b,A\big\}$ of $B$ is $L$--positive.\endslant
%:   Proof of Lemma \LINlem
\Proof Let $C := \lin\big\{b,A\big\}$ and $c \in C$.   Then there exist $\lambda \in \r$ and $a \in A$ such that $c = \lambda b - a$.   If $\lambda = 0$ then, since $0,a \in A \subset \{b\} \cup A$ and $\{b\} \cup A$ is $L$--positive,\break $q_L(c) = q_L(0 - a)  \ge 0$.  If $\lambda \ne 0$ then, since $b \in \{b\} \cup A$, $\lambda^{-1}a \in A \subset \{b\} \cup A$ and $\{b\} \cup A$ is $L$--positive, $q_L(c) = \lambda^2q_L(b - \lambda^{-1}a) \ge 0$.   Thus $q_L \ge 0$ on $C$.   Since $C$ is a linear subspace of $B$, $C$ is $L$--positive.\qed
%:Banach SN duals
\defSection\SNDsec
\medbreak
\centerline{\bf \SNDsec.\quad Banach SN duals}
%:  Definition \DUALdef
\defDefinition \DUALdef
\medskip
\noindent
{\bf Definition \DUALdef.}\enspace Let $(B,L)$ be a Banach SN space and $\big(B^*,\Lt\big)$ also be a Banach SN space.   We say that $\big(B^*,\Lt\big)$ is a \slant Banach SN dual\endslant\ of $(B,L)$ (\cite\LINPOS, Section 2, pp.\ 604--606\endcite) if
$$b \in B\ \and\ b^* \in B^* \qlr \Bra{Lb}{\Lt b^*} = \bra{b}{b^*}.\meqno\DUALone$$
This is equivalent to the definition introduced in \cite\LINPOS, Eq.\ (5), p.\ 605\endcite.
\medskip
In the examples introduced in Example \Hex(a,b,c), the linear map $L$ is orthogonal, from which $(B,L)$ is a Banach SN dual of itself.   However, not all Banach SN spaces have Banach SN duals --- it was shown in \cite\LINPOS, Lemma 8.1, p.\ 614\endcite\ that if $(B,L)$ is a Banach SN space with a Banach SN dual then $L$ is necessarily an isometry, consequently if $B$ is nonzero then $(B,0)$ is a Banach SN space without a Banach SN dual.   We will have a significant example of a Banach SN space with a Banach SN dual in Example \EEex. 
\medskip
We suppose from now on that $(B,L)$ is a Banach SN space with a Banach SN dual $\big(B^*,\Lt\big)$.   We note then from (\DUALone) that,
$$b \in B \qlr \qLt(Lb) = \half\Bra{Lb}{\Lt Lb} = \half\bra{b}{Lb} = q_L(b).\meqno\DUALtwo$$
\par
We will use the following result in Corollary \SURcor.
%:  Lemma \LAlem
\defLemma \LAlem
\medbreak
\noindent
{\bf Lemma \LAlem.}\enspace\slant Let $L(B) = B^*$ and $A$ be a maximally $L$--positive subset of $B$.   Then
$$c^* \in B^* \qlr \infn_{a \in A}\qLt(La - c^*) \le 0.$$\endslant
%:   Proof of Lemma \LAlem
\Proo Let $c^* \in B^*$, and $b \in B$ be chosen so that $Lb = c^*$.   Then, from (\DUALtwo) and (\QMAXone),
$$\infn_{a \in A}\qLt(La - c^*) = \infn_{a \in A}\qLt(La - Lb) = \infn_{a \in A}q_L(a - b) \le 0.\eqno\qed$$\par 
In order to simplify notation in what follows, we will write $\Lambda:= -\Lt$, so $\(B^*,\Lambda\)$ is also a Banach SN space.   We note then from (\DUALone) and (\DUALtwo) that
$$b \in B\ \and\ b^* \in B^* \qlr \bra{b}{b^*} =  -\bra{Lb}{\Lambda b^*},\meqno\DUALthree$$
and
$$b \in B \qlr -q_L(b) = q_\Lambda(Lb).\meqno\DUALfour$$\par
%
%:Section \LINsec: Closed linear subspaces
\defSection \LINsec
\medbreak
\centerline{\bf \LINsec\quad  Closed linear $L$--positive subspaces}
\medskip
\noindent
We will use the following standard notation and results from convex analysis.   If $f\colon\ B \to \rbar$ is proper and convex and $b \in B$ then the \slant subdifferential of $f$ at $b$\endslant, $\partial f(b)$, is the subset of $B^*$ defined by
$$b^* \in \partial f(b) \iff \all\ c \in B,\ f(b) + \bra{c - b}{b^*} \le f(c).$$
\par
We now give a result on the existence of ``almost horizontal'' subtangents:
%:  Lemma \ALHORlem
\defLemma \ALHORlem
\medbreak
\noindent
{\bf Lemma \ALHORlem.}\enspace\slant Let $f\colon\ B \to \rbar$ be proper, convex, lower semicontinuous and bounded below, and $\eta > 0$. Then there exist $b \in \dom\,f$ and $d^* \in \partial f(b)$ such that $\|d^*\| \le \eta$.\endslant
%:   Proof of Lemma \ALHORlem
\Proof This follows from the Br\o ndsted--Rockafellar Theorem \big(see Br\o ndsted--Rockafellar, \cite\BRON, Lemma, pp.\ 608--609\endcite, Z\u alinescu, \cite\ZBOOK, Theorem 3.1.2, p. 161\endcite, or \cite\HBM, Theorem 18.6,\break p.\ 76\endcite\big).\qed\medskip
For the rest of this section, we shall suppose that $A$ is a closed linear $L$--positive subspace of $B$ and $d \in B$.   We define the functions $q_L^{A,d}\colon\ B \to \rbar$ and $q_L^A\colon\ B \to \rbar$ by
$$q_L^{A,d}(b) := \cases{q_L(b)&(if $b \in A - d$);\cr \infty&(otherwise),} \quand q_L^A(b) := q_L^{A,0}(b) = \cases{q_L(b)&(if $b \in A$);\cr \infty&(otherwise).}$$
%:  Lemma \QCSTARlem
\defLemma \QCSTARlem
\par
\noindent
{\bf Lemma \QCSTARlem}\enspace\slant $q_L^{A,d}$ is proper and convex and
$$\partial q_L^{A,d}(b) = \cases{Lb + A^0&(if $b \in A - d$);\cr \emptyset&(otherwise).}$$\endslant
%:   Proof of Lemma \QCSTARlem
\Proo See \cite\LINPOS, Lemmas 5.1 and 5.2, p.\ 609\endcite.\qed
\medskip
The following result is well known:   
%:  Lemma \NORMlem
\defLemma \NORMlem
\medskip
\noindent
{\bf Lemma \NORMlem}\enspace\slant Let $\eta > 0$, $g\colon\ B \to \r$ be defined by $g := 2\eta\|\cdot\|$, $b \in B$ and $b^* \in B^*$.   Then $b^* \in \partial g(b) \iff \|b^*\| \le 2\eta$ and $\bra{b}{b^*} = 2\eta\|b\|$. \endslant
%:  Remark \INFLrem
\defRemark \INFLrem
\medbreak
\noindent
{\bf Remark \INFLrem.}\enspace The technical result that follows is the central result of this paper.   At this time, we point out that the main difference between it and its antecedent, \cite\LINPOS, Theorem 5.3(a), pp.\ 609--610\endcite, is that we use a penalty function of the form $2\eta\|\cdot\|$ rather than $\half\|\cdot\|^2$.   We will explain the implications of this in Remark \INFTrem.
%:  Lemma \INFlem
\defLemma \INFlem
\medbreak
\noindent
{\bf Lemma \INFlem.}\enspace\slant Let $A$ be a closed linear $L$--positive subspace of $B$, $A^0$ be $\Lambda$--positive,\break $d \in B$, $c^* \in B^*$, $\delta := \infn_{a \in A}\|a - d\|$ and $\eta > 0$.   Then there exists $b \in A - d$ such that $\eta\delta \le q_\Lambda(Lb - c^*) +  3\|c^*\|\eta + 5\eta^2$.\endslant 
%:   Proof of Lemma \INFlem 
\Proof Let $g$ be as in Lemma \NORMlem, and write $f := q_L^{A,d} + g - c^*$.   If $\inf_Bf = -\infty$ then we choose $b \in \dom\,f = A - d$ such that $f(b) \le q_\Lambda(c^*) + 3\|c^*\|\eta + 5\eta^2$.   From (\DUALfour) and (\DUALthree),
$$\eqalign{\eta\delta
&\le 2\eta\delta \le g(b) = f(b) - q_L^{A,d}(b) + \bra{b}{c^*} = f(b) - q_L(b) + \bra{b}{c^*}\cr
&= f(b) + q_\Lambda(Lb) - \bra{Lb}{\Lambda c^*} = f(b) + q_\Lambda(Lb - c^*) - q_\Lambda(c^*)\cr
&\le q_\Lambda(Lb - c^*) + 3\|c^*\|\eta + 5\eta^2,}$$
which gives the required result.   Thus we can and will suppose that $\inf_Bf > -\infty$.   From (\QCONTone), $q_L$ is continuous on $B$.   Since $A$ is closed in $B$, $f$ is lower semicontinuous on $B$.   From Lemma \ALHORlem, there exist $b \in \dom\,f = A - d$ and $d^* \in \partial f(b)$ such that $\|d^*\| \le \eta$.   Since $g$ and $c^*$ are continuous, Rockafellar's formula for the subdifferential of a sum \big(see Rockafellar, \cite\FENCHEL, Theorem 3(b), p.\ 85\endcite, Z\u alinescu, \cite\ZBOOK, Theorem 2.8.7(iii), p.\ 127\endcite, or\break \cite\HBM, Theorem 18.1, pp.\ 74--75\endcite\big) provide
$$a^* \in \partial q_L^{A,d}(b)\ \and\ b^* \in \partial g(b)\ \st\ d^* = a^* + b^* - c^*,\ \hbox{thus}\ a^* = c^* + d^* - b^*.\meqno\INFone$$
Lemmas \QCSTARlem\ and \NORMlem\ imply that
$$a^* \in Lb + A^0,\ \|b^*\| \le 2\eta,\ \and\ \bra{b}{b^*} = 2\eta\|b\|.\meqno\INFtwo$$
Using (\INFone), (\INFtwo) and the facts that $\bra{b}{d^*} \le \|b\|\|d^*\| \le \eta\|b\|$ and $\delta \le \|b\|$,
$$\eta\delta + \bra{b}{a^*}
= \eta\delta + \bra{b}{c^*} + \bra{b}{d^*} - \bra{b}{b^*}
\le \eta\delta + \bra{b}{c^*} + \eta\|b\| - 2\eta\|b\| \le \bra{b}{c^*}.$$
Thus, combining with (\DUALthree), 
$$\eta\delta - \bra{Lb}{\Lambda a^*} \le -\bra{Lb}{\Lambda c^*}.\meqno\INFfour$$
From (\INFone) again, $\|a^* - c^*\| \le \|d^*\| + \|b^*\| \le \eta + 2\eta = 3\eta$ thus, applying (\DIFFQ) to $(B^*,q_\Lambda)$,\quad $|q_\Lambda(a^*) - q_\Lambda(c^*)| \le \half\|a^* - c^*\|\|a^* + c^*\| \le \half 3\eta(3\eta + \|2c^*\|) \le 3\|c^*\|\eta + 5\eta^2$,\quad from which
$$q_\Lambda(a^*) \le q_\Lambda(c^*) + 3\|c^*\|\eta + 5\eta^2.\meqno\INFthree$$
From the $\Lambda$--positivity of $A^0$, (\INFtwo), (\INFfour) and (\INFthree),
$$\eqalign{\eta\delta &\le \eta\delta + q_\Lambda(Lb - a^*)
= \eta\delta +  q_\Lambda(Lb) - \bra{Lb}{\Lambda a^*} + q_\Lambda(a^*)\cr
&\le q_\Lambda(Lb) - \bra{Lb}{\Lambda c^*} + q_\Lambda(c^*) + 3\|c^*\|\eta + 5\eta^2
= q_\Lambda(Lb - c^*) + 3\|c^*\|\eta + 5\eta^2.}$$
This completes the proof of Lemma \INFlem.\qed
\medskip
In fact, we will not need Lemma \INFlem\ in full generality.   We will need the following two special cases: Corollary \INFONEcor\ is obtained by setting $c^* = 0$ and noting from (\DUALfour) that $q_\Lambda(Lb) = - q_L(b)$, and Corollary \INFTWOcor\ is obtained by setting $d = 0$ and noting that $\eta\delta \ge 0$.
%:  Corollary \INFONEcor
\defCorollary \INFONEcor
\medbreak
\noindent
{\bf Corollary \INFONEcor.}\enspace\slant Let $A$ be a closed linear $L$--positive subspace of $B$, $A^0$ be $\Lambda$--positive,\break $d \in B$ and $\eta > 0$.   Then there exists $b \in A - d$ such that $\eta\infn_{a \in A}\|a - d\| + q_L(b) \le 5\eta^2$.\endslant 
%:  Corollary \INFTWOcor
\defCorollary \INFTWOcor
\medbreak
\noindent
{\bf Corollary \INFTWOcor.}\enspace\slant Let $A$ be a closed linear $L$--positive subspace of $B$, $A^0$ be $\Lambda$--positive,\break $c^* \in B^*$ and $\eta > 0$.   Then there exists $b \in A$ such that $0 \le q_\Lambda(Lb - c^*) +  3\|c^*\|\eta + 5\eta^2$.\endslant 
%:  Remark \INFTrem
\defRemark \INFTrem
\medbreak
\noindent
{\bf Remark \INFTrem.}\enspace With no extra effort, one can strengthen the conclusion of Lemma \INFlem\ to \slant there exists $b \in A - d$ such that $\eta\|b\| \le q_\Lambda(Lb - c^*) + 3\|c^*\|\eta + 5\eta^2$\endslant.  Theorem \INFthm(a) below was proved in \cite\LINPOS, Theorem 5.3(d), pp.\ 609--611\endcite.   Theorem \INFthm(b) depends on Lemma \POLlem, of which we were not aware when \cite\LINPOS\endcite\ was written.  Theorem \INFthm(c) is a considerable generalization of the critical part of \cite\LINPOS, Corollary 5.4, p.\ 611\endcite, which assumed the additional condition that,\quad for all $d^* \in B^*$, $\infn_{d \in B}\big[\qLt(d^* - Ld) + \half\|d^* - Ld\|^2\big] \le 0$.\quad   We do not know exactly why the change of penalty function mentioned in Remark \INFLrem\ seems to produce this stronger result.   (Of course, $\Lambda$--positivity is equivalent to $\Lt$--negativity.)  
%:  Theorem \INFthm
\defTheorem \INFthm
\medbreak
\noindent
{\bf Theorem \INFthm.}\enspace\slant Let $A$ be a closed linear $L$--positive subspace of $B$ and $A^0$ be $\Lambda$--positive.  Then:
\smallskip
\noindent
%:  Theorem \INFthm(a)
{\rm(a)}\enspace $A$ is maximally $L$--positive.
\smallskip
\noindent
%:  Theorem \INFthm(b)
{\rm(b)}\enspace $A^0$ is maximally $\Lambda$--positive.
\smallskip
\noindent
%:  Theorem \INFthm(c)
{\rm(c)}\enspace For all $c^* \in B^*$, $\inf_{a \in A}\qLt(La - c^*) \le 0$. 
\endslant 
%:   Proof of Theorem \INFthm(a)
\Proof
(a)\enspace Let $d \in B$ and $\{d\} \cup A$ be $L$--positive.  Let $\eta > 0$.   Corollary \INFONEcor\ provides us with $b \in A - d$ such that $\eta\infn_{a \in A}\|a - d\| + q_L(b) \le 5\eta^2$.   By hypothesis, $q_L(b) \ge 0$, hence $\eta\infn_{a \in A}\|a - d\| \le 5\eta^2$.   Thus $\infn_{a \in A}\|a - d\| \le 5\eta$.   Since this holds for all $\eta > 0$ and $A$ is closed, $d \in A$.   This completes the proof that $A$ is maximally $L$--positive.
\smallbreak
%:   Proof of Theorem \INFthm(b)
(b)\enspace Let $c^* \in B$ and $\{c^*\} \cup A^0$ be $\Lambda$--positive.   Write $Z$ for the linear subspace\break $\big\{b \in B\colon\ \bra{b}{c^*} = 0\big\}$ of $B$.   $A \cap Z$ is obviously a closed linear subspace of $B$ and, since $A \cap Z \subset A$, $A \cap Z$ is $L$--positive.   However, from Lemma \POLlem, $(A \cap Z)^0 = \lin\{c^*,A^0\}$ and, from Lemma \LINlem\ \(applied to the SN space $(B^*,\Lambda)$\), $\lin\big\{c^*,A^0\big\}$ is $\Lambda$--positive.    Thus,\break applying (a) with $A$ replaced by $A \cap Z$, $A \cap Z$ is maximally $L$--positive.   Since $A$ is $L$--positive and $A \cap Z \subset A$, in fact $A \cap Z = A$, that is to say $A \subset Z$.   Consequently $\bra{A}{c^*} = \{0\}$, that is to say $c^* \in A^0$.   This completes the proof that $A^0$ is maximally $\Lambda$--positive.
\smallbreak
%:   Proof of Theorem \INFthm(c)
(c)\enspace Let $c^* \in B^*$.  Let $\eta > 0$.   Corollary \INFTWOcor\ provides us with $b \in A$ such that $0 \le q_\Lambda(Lb - c^*) +  3\|c^*\|\eta + 5\eta^2$. Equivalently,  $\qLt(Lb - c^*) \le 3\|c^*\|\eta + 5\eta^2$.   The result now follows by letting $\eta \to 0$.\qed
%:  Remark \INFCrem
\defRemark \INFCrem
\medbreak
\noindent
{\bf Remark \INFCrem.}\enspace Our next result is an improved version of \cite\LINPOS, Corollary 5.4, p.\ 611\endcite.   The critical parts are the implications (\INFCone)$\lr$(\INFCtwo) and (\INFCone)$\lr$(\INFCfour), which are true without any\break additional conditions.   As pointed out in Remark \INFTrem, in \cite\LINPOS\endcite, we assumed the additional condition that,\quad for all $d^* \in B^*$, $\infn_{d \in B}\big[\qLt(d^* - Ld) + \half\|d^* - Ld\|^2\big] \le 0$\quad for the first implication, and we did not establish the second implication. 
%:  Corollary \INFcor
\defCorollary \INFcor
\medbreak
\noindent
{\bf Corollary \INFcor.}\enspace\slant Let $A$ be a closed linear $L$--positive subspace of $B$.  Then the conditions\endslant\ (\INFCone)---(\INFCfour) are equivalent.\slant
$$A^0\ \hbox{is}\ \Lambda\hbox{--positive.}\meqno\INFCone$$
$$A\ \hbox{is maximally}\ L\hbox{--positive}\ \and,\ \all\ c^* \in B^*,\ \infn_{a \in A}\qLt(La - c^*) \le 0.\meqno\INFCtwo$$
$$\All\ c^* \in B^*,\ \infn_{a \in A}\qLt(La - c^*) \le 0.\meqno\INFCthree$$
$$A^0\ \hbox{is maximally}\ \Lambda\hbox{--positive.}\meqno\INFCfour$$
\endslant
%:   Proof of Corollary \INFcor
\Proo It is clear from Theorem \INFthm(a,c) that (\INFCone)$\lr$(\INFCtwo), and it is obvious that (\INFCtwo)$\lr$\break(\INFCthree).  If (\INFCthree) is true then, for all $c^* \in A^0$, (\DUALone) and (\DUALtwo) imply that
$$\infn_{a \in A}q_L(a) + \qLt(c^*) = \infn_{a \in A}\[q_L(a) - \bra{a}{c^*} + \qLt(c^*)\] = \infn_{a \in A}\qLt(La - c^*) \le 0.$$
Since $\infn_{a \in A}q_L(a) \ge 0$, it follows that $\qLt(c^*) \le 0$, from which $q_\Lambda(c^*) \ge 0$, and the linearity of $A^0$ gives (\INFCone).   Thus (\INFCone)---(\INFCthree) are equivalent.   Clearly, (\INFCfour)$\lr$(\INFCone), and the reverse implication is immediate from Theorem \INFthm(b).\qed
%:  Corollary \SURcor
\defCorollary \SURcor
\medbreak
\noindent
{\bf Corollary \SURcor.}\enspace\slant  Let $A$ be a closed linear $L$--positive subspace of $B$ and $L(B) = B^*$.   Then $A$ is maximally $L$--positive if, and only if, $A^0$ is $\Lambda$--positive.\endslant
%:   Proof of Corollary \SURcor
\Proof  If $A^0$ is $\Lambda$--positive then Theorem \INFthm(a) implies directly that $A$ is maximally $L$--positive.   If, conversely, $A$ is maximally $L$--positive then Lemma \LAlem\ implies that  (\INFCthree) is satisfied, and it then follows from Corollary \INFcor\ that $A^0$ is $\Lambda$--positive.  \qed
%:Monotone sets, type (NI) and linear monotone subspaces
\defSection \LINMONsec
\medbreak
\centerline{\bf \LINMONsec.\quad Monotone sets, type (NI) and linear monotone subspaces}
\medskip
\noindent
We suppose in this section that $E$ is a nonzero Banach space.
%:  Example \EEex
\defExample \EEex 
\medbreak
\noindent
{\bf Example \EEex.}\enspace Let $B := E \times E^*$ and, for all $(x,x^*) \in B$, $\|(x,x^*)\| := \sqrt{\|x\|^2 + \|x^*\|^2}$.   We represent $B^*$ by $E^* \times E\dbs$, under the pairing $\Bra{(x,x^*)}{(y^*,y\dbs)} := \bra{x}{y^*} + \bra{x^*}{y\dbs}$, and define $L\colon\ B \to B^*$ by $L(x,x^*) := (x^*,\wh{x})$.   Then $(B,L)$ is a Banach SN space and, for all $(x,x^*) \in B$, $q_L(x,x^*) = \half\[\bra{x}{x^*} + \bra{x}{x^*}\] = \bra{x}{x^*}$.   If $(x,x^*), (y,y^*) \in B$ then we have $q_L\big((x,x^*) - (y,y^*)\big) = q_L(x - y,x^* - y^*) = \bra{x - y}{x^* - y^*}$.   Thus if $A \subset B$ then $A$ is $L$--positive exactly when $A$ is a nonempty monotone subset of $B$ in the usual sense, and $A$ is maximally $L$--positive exactly when $A$ is a maximally monotone subset of $B$ in the usual sense.   We point out that any finite dimensional Banach SN space of the form described here must have \slant even\endslant\ dimension, and that there are many Banach SN spaces  of finite odd dimension with Banach SN duals.   See \cite\SSDMON, Remark 6.7, p.\ 246\endcite.
\par
As usual, the dual norm on $B^* = E^* \times E\dbs$ is given by  $\|(y^*,y\dbs)\| := \sqrt{\|y^*\|^2 + \|y\dbs\|^2}$. By analogy with the analysis above, we define $\Lt\colon\ B^* \to B\dbs$ by $\Lt(y^*,y\dbs) = \big(y\dbs,\wh{y^*}\big)$. Then $\big(B^*,\Lt\big)$ is a Banach SN space, and, for all $(y^*,y\dbs) \in B^*$, $\qLt(y^*,y\dbs) = \bra{y^*}{y\dbs}$.    Now let $b = (x,x^*) \in B$ and $b^* = (y^*,y\dbs) \in B^*$.   Then $\Bra{Lb}{\Lt b^*} = \Bra{(x^*,\wh{x})}{\big(y\dbs,\wh{y^*}\big)} = \bra{x^*}{y\dbs} + \Bra{\wh{x}}{\wh{y^*}} = \bra{x^*}{y\dbs} + \Bra{x}{y^*} = \bra{x}{y^*} + \bra{x^*}{y\dbs} = \Bra{(x,x^*)}{(y^*,y\dbs)} = \bra{b}{b^*}$.   Thus (\DUALone) is satisfied, and so $\big(B^*,\Lt\big)$ is a Banach SN dual of $(B,L)$.  As in the general case, we define $\Lambda(y^*,y\dbs) := -\Lt(y^*,y\dbs) = -\big(y\dbs,\wh{y^*}\big)$. 
\medskip
Definition \NIdef\ first appeared in \cite\RANGE, Definition 10, p.\ 183\endcite.   It was thought at first that this was a weak definition, but it was proved by Marques Alves and Svaiter in \cite\ASD, Theorem 4.4, pp.\ 1084--1085\endcite, that \slant if $A$ is a maximally monotone subset of $E \times E^*$ of type (NI) then $A$ is of ``type (D)''\endslant\ (the opposite implication is obviously true).   This result was extended in \cite\SSDMON, Theorem 9.9(a), pp.\ 254--255\endcite, where it was proved that \slant if $A$ is a maximally monotone subset of $E \times E^*$ of type (NI) then $A$ is of ``dense type'' and ``type (ED)''\endslant.   This has a number of consequences, which are detailed in \cite\SSDMON, Theorem 9.9(b--f), pp.\ 254--255\endcite\ and \cite\SSDMON, Theorem 9.10, pp.\ 255--256\endcite.   Finally, it was proved by Bauschke, Borwein, Wang and Yao in \cite\BBWYFP, Theorem 3.1, pp.\ 1878--1879\endcite\ that \slant if $A$ is maximally monotone subset of $E \times E^*$ then $A$ is of type (NI) if, and only if, $A$ is of ``Fitzpatrick--Phelps'' type\endslant.
\medskip
What distinguishes the definition of ``type (NI)'' from the definitions of the other classes of maximally monotone sets mentioned above is that it can be recast easily in the language of Banach SN spaces, as we will see in Lemma \RECASTlem\ below.  
%:  Definition \NIdef 
\defDefinition \NIdef
\medbreak
\noindent
{\bf Definition \NIdef.}\enspace Let $A \subset E \times E^*$.   We say that $A$ is \slant of type (NI)\endslant\ if,
$$\all\ (y^*,y\dbs) \in E^* \times E\dbs,\quad \infn_{(x,x^*) \in A}\bra{x^* - y^*}{\wh x - y\dbs} \le 0.$$
\par
Here is the promised reformulation of Definition \NIdef\ in the notation of Banach SN spaces.   We use the conventions introduced in Example \EEex.
%:  Lemma \RECASTlem
\defLemma \RECASTlem
\medbreak
\noindent
{\bf Lemma \RECASTlem}\enspace\slant Let $A \subset B := E \times E^*$.   Then $A$ is of type (NI) if, and only if, for all $c^* \in B^*$, $\infn_{a \in A}\qLt(La - c^*) \le 0$.\endslant
%:   Proof of Lemma \RECASTlem
\Proof This is immediate from the formulae for $\qLt$ and $L$ given in Example \EEex.\qed
\medskip
Let $A$ be a linear subspace of $E \times E^*$.  The \slant adjoint subspace\endslant, $A^\T$, of $E\dbs \times E^*$, is defined by:\quad  $(y\dbs,y^*) \in A^\T \iff \hbox{for all}\ (a,a^*) \in A,\ \bra{a}{y^*} = \bra{a^*}{y\dbs}$.\quad This definition goes back at least to Arens in \cite\ARENS\endcite.   (We use the notation ``$A^\T$'' rather than the more usual ``$A^*$'' to avoid confusion with the dual space of $A$.)   It is clear that
$$(y\dbs,y^*) \in A^\T \iff (y^*,-y\dbs) \in A^0.$$  
Now, for all $(y^*,y\dbs) \in B^* = E^* \times E\dbs$, $q_\Lambda(y^*,y\dbs) = -\bra{y^*}{y\dbs}$, from which it follows easily by tracking through the changes of sign and the order of the  variables that $A^\T$ is a monotone subset of $E\dbs \times E^*$ if, and only if, $A^0$ is a $\Lambda$--positive subset of $B^*$, and $A^\T$ is a maximally monotone subset of $E\dbs \times E^*$ if, and only if, $A^0$ is a maximally $\Lambda$--positive subset of $B^*$.
\medskip
We will now see how the results of Section \LINsec\ translate into our present context.   Our first result, which we have not seen in the literature, is about ``automatic maximality'':
%:  Theorem \MONthm
\defTheorem \MONthm
\medbreak
\noindent
{\bf Theorem \MONthm.}\enspace\slant Let $A$ be a closed linear monotone subspace of $E \times E^*$ and $A^\T$ be a monotone subspace of $E\dbs \times E^*$.   Then $A$ and $A^\T$ are maximally monotone in their respective spaces.\endslant
%:   Proof of Theorem \MONthm
\Proof This is immediate from Theorem \INFthm(a,b).\qed 
%:  Remark \MONrem
\defRemark \MONrem
\medbreak
\noindent
{\bf Remark \MONrem.}\enspace The implication (\MONone)$\lr$(\MONtwo) below was shown by Bauschke, Borwein, Wang and Yao in a two--stage process as follows:   It was first shown in \cite\BBWYLIN, Theorem 3.1(iii)$\lr$(ii), pp.\ 6--9\endcite\ that if $A$ is a maximally monotone linear subspace of $E \times E^*$ and (\MONone) is satisfied then $A$ is of type (NI), and it was subsequently shown in \cite\BBWYBB, Propositions 3.1--2, pp. 4955--4959\endcite\ that if $A$ is a closed monotone linear subspace of $E \times E^*$ and (\MONone) is satisfied then $A$ is maximally monotone.   Both of these proofs are algebraically quite complicated, and depend ultimately on the penalty function\quad $(x,x^*) \mapsto \half\|x\|^2 +  \half\|x^*\|^2$.\quad See Remark \INFLrem.
\par
As far as we know, the implication (\MONone)$\lr$(\MONfour) in this generality has not been observed before.   It was proved in the reflexive case in Yao, \cite\YAO, Theorems 18.4--5, p.\ 394--395\endcite.  It was proved for general Banach spaces and ``skew'' subspaces by Bauschke, Borwein, Wang and Yao in \cite\BBWYBB, Corollary 4.4, p. 4961\endcite\ and for subspaces whose domain (i.e., the projection onto $E$) is also closed in \cite\BBWYBB, Theorem 5.5\((v)$\niff$(vi)\), pp. 4965--4968\endcite.
%:  Corollary \MONcor
\defCorollary \MONcor
\medbreak
\noindent
{\bf Corollary \MONcor.}\enspace\slant Let $A$ be a closed linear monotone subspace of $E \times E^*$.  Then the four conditions below are equivalent.
$$A^\T\ \hbox{is a monotone subspace of}\ E\dbs \times E^*.\meqno\MONone$$
$$A\ \hbox{is maximally monotone of type (NI).}\meqno\MONtwo$$
$$A\ \hbox{is of type (NI).}$$
$$A^\T\ \hbox{is a maximally monotone subspace of}\ E\dbs \times E^*.\meqno\MONfour$$
\endslant
%:   Proof of Corollary \MONcor
\Proo This is immediate from Corollary \INFcor\ and Lemma \RECASTlem.\qed
\medbreak
Corollary \BBcor\ below appears in Brezis--Browder \cite\BB, Theorems 1--2, pp.\ 32--39\endcite.   (The result actually proved in \cite\BB, Theorem 1\endcite\ is stronger than that which appears in Corollary \BBcor.)
%:  Corollary \BBcor
\defCorollary \BBcor
\medbreak
\noindent
{\bf Corollary \BBcor.}\enspace\slant Let $E$ be reflexive and $A$ be a closed linear monotone subspace of $E \times E^*$.   Then $A$ is maximally monotone if, and only if, $A^\T$ is monotone.\endslant
%:   Proof of Corollary \BBcor
\Proof This follows from Corollary \SURcor\ and the comments in Example \EEex.\qed
%:References
\bigskip
\centerline{\bf References}
\medskip
\nmbr\ARENS
\item{[\ARENS]} R. Arens, \slant Operational calculus of linear relations\endslant, Pacific J. Math. {\bf 11} (1961) 9--23.
\nmbr\BB
\item{[\BB]} H. Brezis and F. E. Browder, {\sl  Linear maximal monotone operators and singular nonlinear integral equations of Hammerstein,  type}, Nonlinear analysis (collection of papers in honor of Erich H. Rothe),  pp. 31--42. Academic Press, New York, 1978.
\nmbr\BBWYLIN
\item{[\BBWYLIN]} H. Bauschke, J. M. Borwein, X. Wang and L. Yao \slant For maximally monotone linear relations, dense type, negative-infimum type, and Fitzpatrick-Phelps type all coincide with monotonicity of the adjoint\endslant, http://arxiv.org/abs/1103.6239v1, posted March 31, 2011.
\nmbr\BBWYFP
\item{[\BBWYFP]} H. Bauschke, J. M. Borwein, X. Wang and L. Yao \slant Every maximally monotone operator of Fitzpatrick-Phelps type is actually of dense type\endslant, Optim Lett {\bf6} (2012), 1875--1881.\par
DOI 10.1007/s11590-011-0383-2.
\nmbr\BBWYBB
\item{[\BBWYBB]} H. Bauschke, J. M. Borwein, X. Wang and L. Yao \slant The Brezis-Browder Theorem in a general Banach space\endslant, Journal of Functional Analysis {\bf262}(2012), 4948--4971.\par
 http://dx.doi.org/10.1016/j.jfa.2012.03.023.
\nmbr\BRON
\item{[\BRON]} A. Br\o ndsted and R.T. Rockafellar, \slant On the
Subdifferentiability of Convex Functions\endslant, Proc. Amer. Math.
Soc. {\bf 16}(1965), 605--611.
\nmbr\GMS
\item{[\GMS]}Y. Garc\'\i a Ramos, J. E. Mart\'\i nez-Legaz and S. Simons, \slant New results on $q$--positivity\endslant, Positivity {\bf16}(2012), 543--563.   DOI 10.1007/s11117-012-0191-7 
\nmbr\KN
\item{[\KN]}J. L. Kelley, I. Namioka, and co-authors, \slant Linear 
Topological Spaces\endslant, D. Van Nostrand Co., Inc.,
Princeton -- Toronto -- London -- Melbourne (1963).
\nmbr\ASD
\item{[\ASD]}M. Marques Alves and B. F. Svaiter, \slant On Gossez type (D) maximal monotone operators.\endslant, J. of Convex Anal., {\bf 17}(2010), 1077--1088.
\nmbr\FENCHEL
\item{[\FENCHEL]}R. T. Rockafellar, \slant Extension of Fenchel's duality theorem for convex functions\endslant, Duke Math. J. {\bf33}(1966), 81--89.
\nmbr\RANGE
\item{[\RANGE]}S. Simons, \slant The range of a monotone operator\endslant, J. Math. Anal. Appl. {\bf 199}(1996), 176--201.
\nmbr\HBM
\item{[\HBM]}-----, \slant From Hahn--Banach to monotonicity\endslant, 
Lecture Notes in Mathematics, {\bf 1693},\break second edition, (2008), Springer--Verlag.
\nmbr\SSDMON
\item{[\SSDMON]}-----, \slant Banach SSD spaces and classes of monotone sets\endslant, J. of Convex Anal. {\bf 18}(2011), 227--258.
\nmbr\LINPOS
\item{[\LINPOS]}-----, \slant Linear $L$--positive sets and their polar subspaces\endslant, Set-Valued and Variational Anal. {\bf20}(2012), 603--615. DOI: 10.1007/s11228-012-0206-3.
\nmbr\YAO 
\item{[\YAO]} L. Yao, \slant The Brezis--Browder theorem revisited and properties of Fitzpatrick functions of order $n$\endslant, Fixed-point algorithms for inverse problems in science and engineering, Springer Optim. Appl., {\bf49}(2011), 391--402.
\nmbr\ZBOOK
\item{[\ZBOOK]} C. Z\u{a}linescu, \slant Convex analysis in general vector spaces\endslant, (2002), World Scientific.
%
%
%Sends the final list of names to the log file
\Signoff
\bye